\documentclass[oneside]{amsart}
\usepackage{amsfonts}
\usepackage{amssymb}
\usepackage{amsxtra}
\usepackage{amstext}
\usepackage[english]{babel}

\newcommand{\R}{\mathbb{R}}

\newcommand{\E}{\mathbb{E}}

\newcommand{\pp}{\mathbb{P}}

\newcommand{\kR}{\mathcal{R}}

\newcommand{\kH}{\mathcal{H}}

\newcommand{\kL}{\mathcal{L}}
\newcommand{\ka}{\mathfrak{a}}
\newcommand{\kh}{\mathfrak{h}}
\newtheorem {defin} {Definition} [section]

\newtheorem {theo} {Theorem} [section]
\newtheorem {cor} {Corollary} [section]
\newtheorem {rem} {Remark} [section]

\newcommand\la{\lambda}

\newcommand{\lan}{\langle}
\newcommand{\ran}{\rangle}

\title[Bounded harmonic functions for the Heckman--Opdam Laplacian]
      {Bounded harmonic functions for the Heckman--Opdam Laplacian}
\author{Bruno Schapira}
\address{D\'epartement de Math\'ematiques, B\^at. 425, Universit\'e Paris-Sud 11, F-91405 Orsay, cedex, France. }
\email{bruno.schapira@math.u-psud.fr}

\begin{document}

% rajouter reference de Furstenberg pour la frontiere de Poisson de G/K... 

\begin{abstract} We describe the set of bounded harmonic functions for the Heckman--Opdam Laplacian, when the multiplicity function is larger than $1/2$. We prove that this set is a vector space of dimension the cardinality of the Weyl group. We give some consequences in terms of the associated hypergeometric functions. 
\end{abstract}
 
\keywords{Trigonometric Dunkl theory, bounded harmonic function, Poisson boundary, mirror coupling}

\subjclass[2000]{33C67; 60J45; 60J50}

\maketitle

% operateurs d entrelacements...qu est ce que ca donne ? 
% cadre Dunkl ? si k=1 pas de fonctions harmoniques, ok, car \pi est harmonique minimale. Pour les autres k, peut etre
% utiliser que les fonctions harm. bornees sont regulieres dans \ka_+, et operateurs d entrelacements, mais preservent-ils la radialite ? d apres Ancona on peut traiter le cas Dunkl radial avec outils "classiques"....   

\section{Introduction}

In this paper we will consider the operator $\kL$ (called here Heckman--Opdam Laplacian) on $\R^n$ defined, for $f$ a $C^2$ function, by  
\begin{eqnarray}
\label{explicitlaplacian} \kL f(x) &=& \Delta f(x)+ \sum_{\alpha
\in \kR^+}k_\alpha \coth \frac{\lan\alpha,x\ran}{2}\partial_\alpha f(x) \\
 &\quad& \nonumber - \sum_{\alpha \in \kR^+}k_\alpha
\frac{|\alpha|^2}{4\sinh^2 \frac{\lan\alpha,x\ran}{2}} \{f(x)-f(r_\alpha
x)\}.
\end{eqnarray} 
Here $\Delta$ is the usual Euclidean laplacian, $\kR$ is a root system, $\kR^+$ its positive part, the $r_\alpha $'s are the orthogonal reflexions associated to the roots and $k$ is a positive function invariant under the action of the $r_\alpha$'s (see the next section). 
We denote by $W$ the Weyl group, i.e. the finite group generated by the $r_\alpha$'s. We denote by $L$ the restriction of 
$\kL$ to the set of $W$-invariant functions. A simpler formula for $L$ is given by 
\begin{eqnarray}
\label{radiallaplacian2} L f(x) = \Delta f(x)+ \sum_{\alpha
\in \kR^+}k_\alpha \coth \frac{\lan\alpha,x\ran}{2}\partial_\alpha f(x).
\end{eqnarray}
Our main results are the two following: 
 
\begin{theo}
\label{theoradial}
Assume that $k\ge 1/2$. Then the set of bounded $W$-invariant harmonic functions for the Heckman--Opdam Laplacian is exactly the set of constant functions. In other words the Poisson boundary of $L$ is trivial. 
\end{theo} 

\begin{theo} 
\label{theononradial}
Assume that $k\ge 1/2$. Then the set of bounded harmonic functions for the Heckman--Opdam Laplacian is a vector space of dimension $|W|$. In other words the Poisson boundary of $\kL$ is $W$. 
\end{theo} 

In the next section we will give a precise definition for the terminology "harmonic function". We shall also discuss some consequences of our results in terms of the Heckman--Opdam hypergeometric functions, which are particular eigenfunctions of the operator $\kL$.  

\vspace{0.2cm}
\noindent The first result (Theorem \ref{theoradial}) was already known for values of $k$ corresponding to the case of symmetric spaces of the noncompact type $G/K$. The second result (Theorem \ref{theononradial}) is new even for these particular values of $k$, but should be also compared to the situation on symmetric spaces. There, according to the fundamental work of Furstenberg \cite{F} (see also \cite{GJT}), the Poisson boundary of the Laplace--Beltrami operator (but also of a large class of random walks) is $K/M$. But it was already observed that in the Heckman--Opdam (also called trigonometric Dunkl) theory the group $W$ often plays the same role than $K$ or $K/M$. First geometrically, since there is a kind of Cartan decomposition: any $x\in \R^d$ can be uniquely decomposed as $w\cdot x^W$, with $x^W$ the radial part of $x$ (lying in the positive Weyl chamber) and $w\in W$. In representation theory also \cite{O}: briefly if $\kH$ is the graded Hecke algebra generated by $W$ and the Dunkl--Cherednik operators (see next section), then $(\kH,W)$ shares some properties of the Gelfand pair $(G,K)$, like the fact that in any irreducible finite-dimensional $\kH$-module the subspace of $W$-invariant vectors is at most $1$-dimensional. So in some sense Theorem \ref{theononradial} is another manifestation (let say at an analytical or probabilistic level) of the strong analogy between $W$ and $K$.

We should add that the hypothesis $k>0$ is probably sufficient to get the results of Theorem \ref{theoradial} and \ref{theononradial}. Here we restrict us to the case $k\ge 1/2$, because then the stochastic process associated with $L$ (or $\kL$) a.s. never hit the walls (the hyperplanes orthogonal to the roots, which correspond to the singularities of $L$), and we need it to be sure that the coupling we use is well defined.

The paper is organized as follows. In the next section we recall all necessary definitions. 
In section \ref{secradial} we prove Theorem \ref{theoradial}, by using the probabilistic technique of mirror coupling. In section \ref{secnonradial} we prove Theorem \ref{theononradial}, by extending the coupling to the non-radial process. Our main tool for this is the skew-product representation from Chybiryakov \cite{Chy}, that we have to adapt to our setting. 

\vspace{0.2cm}
\noindent \textit{Acknowledgments: I warmly thank Marc Arnaudon for having explained to me the technique of mirror coupling, 
and Alano Ancona for enlightening discussions about the regularity of harmonic functions. }

\section{Preliminaries}
Let $\ka$ be a Euclidean vector space of dimension $n$, equipped
with an inner product $\lan\cdot ,\cdot \ran$, and denote by $\kh:=\ka +i\ka$ its complexification. 
We consider $\kR \subset \ka$ an integral
root system (see \cite{Bou}). We choose a subset of positive roots $\kR^+$. 
Let $\alpha^\vee=2\alpha/|\alpha|^2$
be the coroot associated to a root $\alpha$ and let
$$r_\alpha(x)=x-\lan\alpha^\vee,x\ran\alpha,$$
be the corresponding orthogonal reflection. Remember that $W$ denotes the
Weyl group associated to $\kR$, i.e. the group generated by the
$r_\alpha$'s. 
Let $k\ :\ \kR \rightarrow [1/2,+\infty)$ be a multiplicity
function, which by definition is $W$-invariant. We set 
$$\rho=\frac{1}{2}\sum_{\alpha \in \kR^+}k_\alpha \alpha.$$
Let
$$\ka_+ = \{x \mid \forall \alpha \in \kR^+,\ \lan\alpha,x\ran>0\},$$
be the positive Weyl chamber. Let also $\overline{\ka_+}$ be its
closure, $\partial \ka_+$ its boundary and
$\ka_{\text{reg}}$ the subset of regular elements in $\ka$,
i.e. those elements which belong to no hyperplane $\{\alpha=0\}$. 
As recalled in the introduction any $x\in \ka$ can be uniquely decomposed as $x=w x^W$, with $x^W\in \overline{\ka_+}$ and $w\in W$. 
We call $x^W$ the radial part of $x$ and $w$ its angular part.

For $\xi \in \ka$, let $T_\xi$ be the Dunkl--Cherednik
operator \cite{C}. It is defined, for $f\in C^1(\ka)$, and $x\in
\ka_{\text{reg}}$, by
$$T_\xi f(x)=\partial_\xi
f(x) + \sum_{\alpha \in \kR^+}k_\alpha
\frac{\lan\alpha,\xi\ran}{1-e^{-\lan\alpha,x\ran}}\{f(x)-f(r_\alpha x)
\}-\lan\rho,\xi\ran f(x).$$ 
The Dunkl-Cherednik operators form a commutative family of
differential-difference operators (see \cite{C} or \cite{O}). The
Heckman--Opdam Laplacian $\kL$ is also given by the formula 
$$\kL+|\rho|^2=\sum_{i=1}^{n} T_{\xi_i}^2,$$
where $\{\xi_1,\dots,\xi_n\}$ is any orthonormal basis of $\ka$.

Let $\la \in \kh$. We denote by $F_\la$ the unique (see \cite{HO}, \cite{O}) analytic
$W$-invariant function on $\ka$, which satisfies the differential
equations 
$$p(T_\xi)F_\la=p(\la)F_\la \text{ for all W-invariant polynomials }p$$
and which is normalized by $F_\lambda(0)=1$ (in particular $\kL
F_\la=(\lan\la,\la\ran-|\rho|^2) F_\la$). We denote by $G_\lambda$ the unique
analytic function on $\ka$, which satisfies the
differential-difference equations (see \cite{O})
\begin{eqnarray}
\label{equations} T_\xi G_\la = \lan\la,\xi\ran G_\la \text{ for all }\xi
\in \ka,
\end{eqnarray}
and which is normalized by $G_\lambda(0)=1$. These functions are related by the formula: 
\begin{eqnarray}
\label{FG}
F_\la(x)=\frac{1}{|W|} \sum_{w\in W} G_\la(wx),
\end{eqnarray}
for all $x\in \ka$ and all $\la \in \kh$.

It was shown in \cite{Sch2} that $\frac{1}{2}\kL$ and $\frac{1}{2}L$ are generators of Feller semi-groups that we shall 
denote respectively by $(P_t,t\ge 0)$ and $(P^W_t,t\ge 0)$. We will use the following definition for harmonic functions:   

\begin{defin} 
\label{defharmonic}
A bounded or nonnegative function $h:\ka \to \R$ is called harmonic if it is measurable and satisfies $P_th=h$ for all $t>0$.  
\end{defin}

\begin{rem} \emph{It is well known that if $h$ is a $C^2$ function such that $\kL h=0$, then $h$ is harmonic in the sense of Definition \ref{defharmonic}. 
Inversely Corollary \ref{coroG} below shows, when $k\ge 1/2$, that any bounded harmonic function is regular, thus satisfies $\kL h=0$. On the other hand, 
it is a general fact (which applies for any $k>0$), that bounded $W$-invariant harmonic functions are regular in $\ka_+$, but we will not use this fact here.}
\end{rem} 

Observe that by definition $F_\rho$ is a $W$-invariant harmonic function. Moreover it is known (see \cite{Sch2} Remark 3.1) that it is bounded. So Theorem \ref{theoradial} shows that in fact $F_\rho$ is constant equal to $1$. Similarly the functions $G_{w\rho}$'s, for $w\in W$, are harmonic and also bounded. This last property follows from Formula \eqref{FG}, since the $G_{w\rho}$'s are real positive (see \cite{Sch2} Lemma 3.1). In fact one has the following 
\begin{cor}
\label{coroG} 
If $k\ge 1/2$, then any bounded harmonic function is a linear combination of the $G_{w\rho}$'s, $w\in W$. 
\end{cor} 
\begin{proof}
The only thing to prove is that the $G_{w\rho}$'s are linearly independent. This results from the fact that they are all eigenfunctions of the Dunkl--Cherednik operators but for different eigenvalues. More precisely, assume that for some real numbers $(c_w)_{w\in W}$, we have 
$$\sum_{w\in W} c_w G_{w\rho}=0.$$
By applying then the operators $p(T_\xi)$, with $p$ polynomial, we get 
$$\sum_{w\in W} c_w p(w\rho)G_{w\rho}=0 \quad \textrm{for all p}.$$
>From this, and the fact that $G_{w\rho}(0)=1$ for all $w$, it is easily seen that we must have $c_w=0$ for all $w$.
\end{proof}

\section{The $W$-invariant case: proof of Theorem \ref{theoradial}}
\label{secradial}
In this section we shall prove Theorem \ref{theoradial}. For this we will use the stochastic process $(X^W_t,t\ge 0)$ associated with $L$, called radial HO-process, and the so-called mirror coupling technique.

First it is known \cite{Sch1} that $X^W$ is a strong solution of the SDE: 
$$     
X^W_t=x+B_t + V^1_t
$$
where $(B_t,t\ge 0)$ is a Brownian motion on $\ka$ and 
$$
V^1_t:=\sum_{\alpha\in \kR^+} k_\alpha \alpha \int_0^t \coth \lan \alpha, X^W_s\ran \ d s.
$$
Moreover when $k\ge 1/2$, $X^W$ a.s. takes values in $\ka_+$, or in other words it never reaches $\partial \ka_+$ (see \cite{Sch1}).  
Now if $x,y\in \ka_+$, we define the couple $((X^W_t,Y^W_t),t\ge 0)$ as follows. Set $T=\inf\{s \mid X^W_s=Y^W_s\}$. Then by definition $X^W$ is as above, and $(X^W,Y^W)$ is the unique solution of the SDE:  
\begin{eqnarray} 
\label{SDE}
(X^W_t,Y^W_t)=(x,y) + (B_t,B'_t)  + (V^1_t,V^2_t), \quad \textrm{for } t<T,
\end{eqnarray}
where $dB'_t=r_t dB_t$, with $r_t$ the orthogonal reflexion with respect to the hyperplane orthogonal to the vector $Y^W_t-X^W_t$ (in particular Levy criterion shows that $B'$ is a Brownian motion), and 
$$
V^2_t:=\sum_{\alpha\in \kR^+} k_\alpha \alpha \int_0^t \coth \lan \alpha, Y^W_s\ran \ d s.
$$
For $t\ge T$, we set $Y^W_t=X^W_t$. The existence of this coupling is guaranteed by the fact that the SDE \eqref{SDE} has locally regular coefficients. We define also $Z^W$ by  
$$Z^W_t:=Y^W_t-X^W_t,$$
and set $z^W_t= |Z^W_t|$. It is known \cite{Sch1} that a.s. $X^W_t/t \to \rho$, and thus that $\lan\alpha,X^W_t\ran \sim \lan\rho,\alpha\ran t$, for all $\alpha \in \kR^+$. From this we see that a.s.  $\sup_{t\ge 0}|V^2_t-V^1_t|<+\infty$. Then Tanaka formula (\cite{RY} p.222) shows that 
$$
z^W_t =\gamma_t + v_t, \quad \textrm{for } t<T,
$$
with $\gamma$ a one-dimensional Brownian motion and a.s. $\sup_{t\ge 0}|v_t|<+\infty$. In particular $T$ is a.s. finite.

The end of the proof is routine now. Assume that $h$ is a bounded $W$-invariant harmonic function. Then it is well known, and not difficult to show, that $(h(X^W_t),t\ge 0)$ as well as $(h(Y^W_t),t\ge 0)$ are bounded martingales. Thus they are a.s. converging toward some limiting (random) values, respectively $l$ and $l'$. Since a.s. $X^W_t=Y^W_t$ for $t$ large enough, we have a.s. $l=l'$. 
Then usual properties of bounded martingales show that 
$$
h(x) = \E[l] = \E[l']=h(y).
$$
Since this holds for any $x,y \in \ka_+$, this proves well that $h$ is constant.  \hfill $\square$

\section{The non $W$-invariant case: proof of Theorem \ref{theononradial}}
\label{secnonradial}
In order to prove Theorem \ref{theononradial}, the first idea is to extend the previous coupling to the full process $(X_t,t\ge 0)$ with 
semi-group $(P_t,t\ge 0)$. For this our tool will be the skew-product representation founded by Chybiryakov \cite{Chy} (see \cite{GaY} and \cite{Chy2} for the one-dimensional case). Actually Chybiryakov dealt with Dunkl processes, so we shall first mention the changes needed to adapt his proof to the present setting, and then explain how to combine this representation with the coupling from the previous section. 

\subsection{Skew-product representation and extension of the coupling} 
\label{sprod}
The skew-product representation gives a constructive way to define $X$ starting from $X^W$, by adding successively jumps in the direction of the roots. 
Let us sketch the main steps of the construction (for more details see \cite{Chy}). First one fixes arbitrarily an order for the positive roots: $\alpha_1,\dots,\alpha_{|\kR^+|}$. Then for each $j \in [1,|\kR^+|]$, set 
$$
\kL^jf(f):= Lf(x) - \sum_{i\le j}c_{\alpha_i}(x) \{f(x)-f(r_{\alpha_i}
x)\},
$$
where for any root $\alpha$, 
$$
c_\alpha(x):= k_{\alpha}
\frac{|\alpha|^2}{4\sinh^2 \frac{\lan\alpha,x\ran}{2}}.
$$
Decide also that $\kL^0=L$. 
Set 
$$
\widetilde{\kL}^j f(x) : = c_{\alpha_j}^{-1}(x) \kL^j f(x),
$$
and 
$$
\kL^{j,j+1}f(x) :=  c_{\alpha_{j+1}}^{-1}(x) \kL^j f(x).
$$
The goal is to define inductively a sequence of processes $(X^j(t),t\ge 0)$, $j=0,\dots,|\kR^+|$, associated to the operators $\kL^j$'s. 
First $X^0$ is just the radial HO-process considered in the previous section. Next assume that $\kL^j$ is the generator of a Markov process $(X^j(t),t\ge 0)$. Then set 
$$
A_t^j=\int_0^t c_{\alpha_{j+1}}(X^j_s)\ ds,
$$
and 
$$
\tau_t^j= \inf\{s\ge 0\mid A_s^j>t\}.
$$
Using the martingale problem characterization one can see that the radial part of $X^j$ is a radial HO-process. Thus for all $\alpha\in \kR^+$, $|\lan \alpha,X^j_t\ran|\ge c t$, for $t$ large enough and $c>0$ some constant.
In particular the increasing process $A^j$ is bounded. Set $T^j=\lim_{t\to +\infty} A^j_t$. Then observe that $\tau_t^j=+\infty$, when $t\ge T^j$. This is essentially the only difference with the Dunkl case considered in \cite{Chy} (where $A^j$ was not bounded and $\tau_t^j$ finite for all $t$). But one can still see 
that if 
$$X^{j,j+1}(t):= X^j(\tau_t^j) \quad t< T^j,$$
then $X^{j,j+1}$, killed at time $T^j$, is solution of the martingale problem associated with $\kL^{j,j+1}$ (see for instance \cite{EK} exercise 15 p.263 and section 6 p.306). The next step is to add jumps to $X^{j,j+1}$ in the direction of the root $\alpha_{j+1}$. Namely one define a new process $\widetilde{X}^j$, also denoted by $X^{j,j+1}*_{\alpha_{j+1}} N$ in \cite{Chy} section 2.5, which is solution of the martingale problem associated with $\widetilde{\kL}^{j+1}$. Roughly $\widetilde{X}^j$ is constructed by gluing several paths, all with law $X^{j,j+1}$ or $r_{\alpha_{j+1}}X^{j,j+1}$, such that for any two consecutive path the starting point of the second is the image of the end point of the first path by the reflexion $r_{\alpha_{j+1}}$. The lengths of the paths are determined by independent exponentially distributed random variables. Here the only minor change is that $\widetilde{X}^j$ explodes at some time, let say $\widetilde{T}^j$. A change of variables shows that 
$$
\lim_{t\to \widetilde{T}^j} \int_0^t c_{\alpha_{j+1}}^{-1}(\widetilde{X}^j(s))\ ds=+\infty.
$$
So for any $t\ge 0$, one can define $\widetilde{A}^j(t)$ as solution of the equation 
$$
t = \int_0^{\widetilde{A}^j(t)}c_{\alpha_{j+1}}^{-1}(\widetilde{X}^j(s))\ ds.
$$
Differentiating this equation one get 
$$
\frac{d}{dt} \widetilde{A}^j(t) = c_{\alpha_{j+1}}(\widetilde{X}^j(\widetilde{A}^j(t))). 
$$ 
Then set $X^{j+1}(t)=\widetilde{X}^j(\widetilde{A}^j(t))$, for all $t\ge 0$. The preceding equation gives  
$$
\widetilde{A}^j(t)= \int_0^t c_{\alpha_{j+1}}(X^{j+1}(s))\ ds,
$$ 
which in turn shows that $X^{j+1}$ is solution of the martingale problem associated with $\kL^{j+1}$, as wanted.

The point now is to combine this construction of $X=X^{|\R^+|}$ with the coupling of the radial process from section \ref{secradial}. 
We first take $(X^0,Y^0)$ with law given by this coupling. Then we define the sequence $((X^j(t),Y^j(t)),t\ge 0)$, 
$j= 1,\dots,|\kR^+|$, simply by following the previous construction for the two coordinates. Actually this coupling is interesting 
only when $X=X^{|\kR^+|}$ and $Y=Y^{|\kR^+|}$ never jump, but this is precisely what we need. Indeed in this case we have $X_t=X^0(t)$ and $Y_t=Y^0(t)$, for all $t\ge 0$, so they coincide a.s. after some finite time.

%%%%%%%%%%%%%%%%%%%%%%%%%%%%%%%%%%%%%%%%%%%%%%%%%%%%%%%%%%%%%%%%%%%%%%%%%%%%%%%%%%%%%%%%%%%%%%%%%%%
\subsection{End of the proof}

For any $x\in \ka$, we denote by $\pp_x$ the law of $(X_t,t\ge 0)$ starting from $x$. 
For $\epsilon \in (0,1)$, set 
$$A_\epsilon := \{z\in \ka \mid \pp_z[X \textrm{ never jumps}] \ge 1-\epsilon\}.$$
We know that the process $(X_t,t\ge 0)$ can jump, so a priori $A_\epsilon \subsetneq \ka$. But we know also \cite{Sch1} that a.s. $X$ eventually stops to jump after some finite random time. This implies that
\begin{eqnarray}
\label{eqjumps}
\lim_{t\to +\infty} \pp_x[X \textrm{ never jumps after time } t]=1,
\end{eqnarray}
for all $x\in \ka$.
But by using the Markov property, we have for all $t>0$,
\begin{eqnarray}
\label{eqjump2}
\nonumber \pp_x[X \textrm{ never jumps after time } t] &=& \E_x\left[ \pp_{X_t}[X \textrm{ never jumps}] \right]\\
                                             &=& \int_\ka \pp_z[X \textrm{ never jumps}]\ d\mu^x_t(z),
\end{eqnarray}
where $\mu^x_t$ is the law of $X_t$ under $\pp_x$. So \eqref{eqjumps} and \eqref{eqjump2} imply that for all $x\in \ka$,
$\mu_t^x(A_\epsilon) \to 1$, when $t\to +\infty$. In particular $A_\epsilon$ is nonempty.  
Moreover, by invariance of $\kL$ under $W$, we know that for any $w\in W$, the law of $(wX_t,t\ge 0)$ under $\pp_x$ is $\pp_{wx}$. 
In particular, for any $w\in W$ and any $\epsilon\in (0,1)$, we have $w(A_\epsilon\cap \ka_+)=A_\epsilon \cap w\ka_+$. Thus all these subsets of $A_\epsilon$ are nonempty as well.

Let now $h$ be some harmonic function. Fix $w\in W$, and take $x,y \in A_\epsilon \cap w\ka_+$. Consider the coupling $((X_t,Y_t),t\ge 0)$ as defined above. Since $(h(X_t),t\ge 0)$ and $(h(Y_t),t\ge 0)$ are bounded martingales, they converge a.s. toward some limits, respectively $l$ and $l'$. 
We already saw that $X^W$ and $Y^W$ a.s. coincide after some time. So if both processes $X$ and $Y$ never jump, they must also coincide after some time, and in this case we have $l=l'$. Since $x,y \in A_\epsilon$, this shows that 
$$|h(x)-h(y)|=|\E[l]-\E[l']| \le 2C\epsilon,$$
where $C=\sup h$. In particular, by completeness of $\R$, for any sequence $(x_\epsilon)_{\epsilon \in (0,1)}$, such that $x_\epsilon \in A_\epsilon\cap w\ka_+$ for all $\epsilon \in (0,1)$, the limit of $h(x_\epsilon)$ when $\epsilon $ tends to $0$ exists, and is independent of the chosen sequence. Call $l_w$ this limit.

For all $t\ge 0$, we denote by $w_t$ the angular part of $X_t$. Since $X$ eventually stops to jump, $(w_t,t\ge 0)$ a.s. converges, i.e. becomes stationary. Then for any $w\in W$, define the function $h_w$ on $\ka$ by 
$$h_w(x)=\pp_x\left[\lim_{t\to +\infty}w_t=w\right].$$
By standard properties of Markov processes, we know that these functions are measurable, and actually it is not difficult to see that they are harmonic. Moreover the above convergence result for harmonic functions shows that these functions 
$h_w$, $w\in W$, are linearly independent. Then set 
$$\tilde{h}(x):= \sum_{w\in W} l_w h_w(x),$$
for all $x\in \ka$. All that remains to do now is to prove that $\tilde{h}=h$. Indeed if this was true, this would prove that the vector space of bounded harmonic function has dimension $|W|$ as wanted. 
By using the martingale property, we have for any $t>0$ 
\begin{eqnarray}
\label{hhtilde}
|h(x)-\tilde{h}(x)|= |\E_x[h(X_t)-\tilde{h}(X_t)]| \le \int_\ka |h(z)-\tilde{h}(z)|\ d\mu_t^x(z).
\end{eqnarray}
We have seen that for all $\epsilon \in (0,1)$, 
\begin{eqnarray}
\label{Aepsilon}
\mu_t^x(A_\epsilon) \to 1
\end{eqnarray}
when $t\to +\infty$. But it is not difficult to see (by using the definition of the $l_w$'s), that for any $\epsilon'>0$, there exists $\epsilon>0$ such that 
$$|h(z)-\tilde{h}(z)|\le \epsilon' \quad \forall z\in A_\epsilon. $$   
Since this holds for any $\epsilon'>0$, \eqref{hhtilde} and \eqref{Aepsilon} show that $h=\tilde{h}$ as wanted. \hfill $\square$

\begin{rem} \emph{We have seen in the previous proof that the family $(h_w)_{w\in W}$ is a basis of the space of bounded harmonic functions. Since the family $(G_{w\rho})_{w\in W}$ is another basis, it would be interesting to know the coefficients relating these two basis. }
\end{rem}

\end{document}